\documentclass[11pt, reqno,twoside]{article}

\usepackage{amsmath,amsfonts,amssymb,amsthm,bbm}
\usepackage{graphics,epsfig}
\usepackage{color}
\usepackage{array}

\newtheorem{theorem}{Theorem}[section]
\newtheorem{proposition}{Proposition}[section]
\newtheorem{lemma}{Lemma}[section]
\newtheorem{corollary}{Corollary}[section]

\numberwithin{equation}{section}

\allowdisplaybreaks

\def\beqlb{\begin{eqnarray}}
\def\eeqlb{\end{eqnarray}}
\def\beqnn{\begin{eqnarray*}}
\def\eeqnn{\end{eqnarray*}}

\def\qed{\hfill$\Box$\medskip}

\begin{document}

\bigskip

\bigskip

\bigskip

\title{{\bf Explicit stationary distribution of the $(L,1)$-reflecting random walk on the half line}\thanks{The project is partially supported by
{the National Natural Science Foundation of China (Grant No. 11131003) and by the Natural Sciences and Engineering Research Council of Canada (Grant No. 315660).}}}

\author{Wenming Hong\footnote{School of Mathematical Sciences
\& Laboratory of Mathematics and Complex Systems, Beijing Normal
University, Beijing 100875, P.R. China. Email: wmhong@bnu.edu.cn}, \ \ \
 Ke Zhou\footnote{School of Mathematical Sciences
\& Laboratory of Mathematics and Complex Systems, Beijing Normal
University, Beijing 100875, P.R. China. Email:zhouke@mail.bnu.edu.cn} \ \ \ and  \ \ \ {Yiqiang Q. Zhao}\footnote{School of Mathematics and Statistics,
Carleton University, Ottawa, Ontario, Canada K1S 5B6. Email: zhao@math.carleton.ca}}
\date{}

\maketitle

\bigskip\bigskip

{\narrower{\narrower

\centerline{\bf Abstract}

\bigskip
In this paper, we consider the $(L,1)$ state-dependent reflecting random walk (RW) on the half line, which is a RW allowing jumps to the left at a maxial size $L$.
For this model, we provide an explicit criterion for (positive) recurrence and an explicit expression for the stationary distribution.
As an application, we prove the geometric tail asymptotic behavior of the stationary distribution under certain conditions. The main tool employed in the paper is the intrinsic branching structure within the $(L,1)$-random walk.

\bigskip


\bigskip


\par}\par}

\bigskip\bigskip


\section{ Introduction and main results\label{s1}}

\subsection{The background and motivation}
We consider the $(L,1)$-reflecting random walk on the half line,  i.e.,  a Markov chain $ \{X_n\}_{n\geq 0}$ on $\mathbb{Z^{+}}=\{0, 1, 2, \ldots \}$ with $X_0=0$ and the transition probabilities $P_{ij}$ specified by
: $P_{01}=p(0)=1$;  for $0< i< L$,
\beqnn
    P_{ij} = \left \{ \begin{array}{ll}
p(i), & \mbox{for $j=i+1$}, \\
q_{i-j}(i), & \mbox{for $0<j<i$}, \\
 \sum_{ k=i}^{L}q_k(i), & \mbox{for $j=0$}, \\
0, & \mbox{otherwise,}
\end{array} \right. \eeqnn
and for $i\geq L$,
\beqnn
P_{ij} = \left \{ \begin{array}{ll}
p(i), & \mbox{for $j=i+1$}, \\
q_{i-j}(i), & \mbox{for $i-L\leq j<i$}, \\
0, & \mbox{otherwise,}
\end{array} \right.
   \eeqnn where $q_{1}(i)+q_{2}(i)+\cdots +q_{L}(i)+p(i)=1$, and $q_{1}(i), q_{2}(i), \cdots, q_{L}(i)\geq 0, 0<p(i)<1$. Obviously, this Markov chain  is irreducible.

For example, when $L=2$, the transition matrix $P$ is given by\beqnn
 P=\left(
             \begin{array}{ccccccc}
               0 & 1 & \\
               q_{1}(1)+q_{2}(1) & 0 & p(1) & \\
               q_{2}(2) & q_{1}(2) & 0 & p(2) & \\
               & q_{2}(3) & q_{1}(3) & 0 & p(3) & \\
               & & \ddots & \ddots & \ddots & \ddots & \\
             \end{array}
           \right).
\eeqnn
in which all unspecified entries are zero.
\

It is well-known that for the $(1,1)$-RW,  the
criteria for the (positive) recurrence and the  expression for the
stationary distribution have been given explicitly, for example
see \cite{K-M:57} and \cite{Kar}). However, for $L >1$, no such
explicit expressions are found to our best knowledge. The aim of
the present paper is to give explicit criteria of the (positive)
recurrence and explicit expressions of  the stationary
distribution for the $(L,1)$-RW. Our method is probabilistic by
using the intrinsic branching structures hidden in the $(L,1)$-RW
(\cite{hw} and \cite{hzh}).  The results obtained in this paper
can be applied to various state-dependent queueing systems.

\
\subsection{Main results}

 First, we recall two results on the (positive)
recurrence and the stationary distribution of a general Markov
chain $X_n$ (e.g., \cite{Kar} or \cite{Durr}). For this purpose,
let $N(i)=\sum_{n=1}^{\infty} 1_{\{X_n=i\}}$ be the number of
visits to state $i$ by the chain, $T_{i}=\inf\{n>0: X_{n}=i\}$ is
the first time for the chain to be in state $i$, and $E^i$ is the
expected value when the walk starts at $X_0=i$.
\begin{description}
\item[Fact~1] $\{X_{n}\}_{n\geq 0}$  is recurrent
 $\Longleftrightarrow E^0N(0)=\infty$ $\Longleftrightarrow
\sum_{j=0}^{\infty}P_{ij}y_{j}=y_{i}$ ($i\geq 0$) have no bounded
nonconstant solution; and

\item[Fact~2] $\{X_{n}\}_{n\geq 0}$  is positive recurrent
$\Longleftrightarrow E^0T_0<\infty$. In this case, the stationary
distribution $\pi(i)$ is given by $\pi(i)=\frac{1}{E^{i}T_{i}}$.
\end{description}
 The idea of the present
paper is to express $E^0N(0)$ and $E^{i}T_{i}$ explicitly through
using the intrinsic branching structure hidden in the $(L,1)$-RW.

Three results are obtained: the first one is a characterization
for recurrence; the second is an expression for the stationary
probability distribution; and the last one is a criterion for the
stationary probabilities to have a geometric decay. For the first
two results, we present them for $L=2$ since the notation for a
general $L$ is very demanding. The last result is presented for a
general $L$.

\subsubsection{\it  Criteria for the recurrence} ~
Let $e_1=(1,0)$, $e_2=(0,1)$, $u=e_1'+e_2'$, the sum of the
transposes of $e_1$ and $e_2$,
 \beqnn
 M_{1}= \left(\begin{array}{cc}
\frac{q_{1}(1)+q_{2}(1)}{p(1)} & 0\\
\frac{1}{p(1)} & 0
\end{array} \right), \quad
    M_{i}= \left(\begin{array}{cc}
\frac{q_{1}(i)}{p(i)} & \frac{q_{2}(i)}{p(i)}\\
1+\frac{q_{1}(i)}{p(i)} & \frac{q_{2}(i)}{p(i)}
\end{array} \right), \;i>1.
\eeqnn
\begin{theorem}\label{t1.1}
Let
\beqnn
\kappa:=\sum_{k=1}^{\infty}e_{1}M_{k}M_{k-1}\ldots M_{1}u,\eeqnn
if $\kappa= \infty$, then the walk $\{X_{n}\}_{n\geq 0}$ is recurrence.
\end{theorem}
\noindent {\bf Remark} Actually, $\kappa$ in the theorem is the expectation number of the visiting times by the random walk at position $0$, i.e., $\kappa=E^0N(0)$ , which can be calculated by the means of the  intrinsic branching structure within the walk and the $M_i$  is the offspring mean matrix.

\subsubsection{\it Criteria for the positive recurrence and stationary distribution}~
Next, we will give the criteria for the positive recurrence and  explicit formula of the stationary distribution based on the Fact 2. Let $P^{i}$ denote the probability when the walk start at $X_0=i$. Define
\beqlb\label{f1.3}
\tau_0=0,~~~~\tau_i=\inf\{n> 0:X_n< i\}, \quad i\ge 1 \eeqlb

\begin{theorem}\label{t1.2} Assume $E^{1}\tau_{1}<\infty $, then the
walk $\{X_{n}\}_{n\geq 0}$ is positive recurrence. Furthermore the
stationary distribution is given by, for $i\geq 0$,
 \beqnn
    \pi(i)=\frac{1}{E^{i}T_{i}}=\frac{1}{p(i)E^{i+1}
\tau_{i+1}+(q_{1}(i)+q_{2}(i)+p(i)P^{i+1}[(i+1,+\infty),i-1])E^{i-1}T_{i}+q_{2}(i)E^{i-2}T_{i-1}+1},
 \eeqnn
where explicit expressions for
$P^{i+1}[(i+1,+\infty),i-1]$, $E^{i}\tau_{i}$ and $E^{i}T_{i+1}$
are given in (\ref{f1.4}), (\ref{f1.6}) and (\ref{f1.7}),
respectively.
\end{theorem}

\subsubsection{\it Tail behavior of the stationary distribution}~
It is well know { that when the transition probabilities
are} state-independent, the tail of the stationary distribution of
the walk has a geometric decay. Here with the
 help of the explicit expression of the stationary
distribution in given in the previous theorem, we can
consider the tail behavior for the state-dependent
case. We notice that the expression of the stationary
distribution in Theorem \ref{t1.2} is given in terms
of the decomposition of the trajectory, from which we
find out that the dominant contribution to the tail asymptotic
behavior of $\pi (i)$ is from $E^{i-1}T_{i}$ ( because of
$E^{1}\tau_{1}<\infty$). With this observation,
characterize the tail behavior of the state-dependent $(L,1)$-RW
as follows.

Let $D:=\{(p,q_{1},q_{2},\cdots q_{L}): p+\sum_{j=1}^{L}q_{j}=1; ~~ \sum_{j=1}^{L}jq_{j}>p;~~\forall j, ~q_{j}\geq 0, p>0 \}.$ $\rho(M)$ is the spectral radius of $M$, and $\lambda_{M}$ is the maximum eigenvalue of it, where
\beqlb\label{m}
M =\left(\begin{array}{ccccc}
\frac{q_{1} }{p } & \frac{q_{2} }{p } & \ldots & \frac{q_{L-1} }{p } & \frac{q_{L} }{p  }\\
1+\frac{q_{1} }{p } & \frac{q_{2} }{p(i)} & \ldots & \frac{q_{L-1} }{p } & \frac{q_{L} }{p }\\
\vdots & \vdots & \ddots & \vdots & \vdots\\
\frac{q_{1} }{p } & \frac{q_{2} }{p } &\cdots & \frac{q_{L-1} }{p } & \frac{q_{L} }{p }\\
\frac{q_{1} }{p } & \frac{q_{2} }{p } & \cdots & 1+\frac{q_{L-1} }{p } & \frac{q_{L} }{p }
\end{array} \right)_{L\times L}.
\eeqlb

\begin{theorem}\label{t1.3}
(1) If $P\in D  $, we have $\lambda_{M}=\rho(M)>1.$  \\
(2) If $P(i):=(p(i),q_{1}(i),q_{2}(i), \cdots, q_{L}(i) )\rightarrow P$, as $i\to\infty$, and $P\in D  $, then the stationary distribution exist; and
$$\lim_{i\rightarrow \infty}\frac{\log\pi({i})}{i}= -\log{\lambda_{M}}.$$
\end{theorem}
\noindent {\bf Remark}~ Here    $D$ is the positive recurrence district  for the state-independent $(L, 1)$-RW, and $M$ is the correspond offspring mean matrix. The theorem says that if the transition probability of the state-dependent $(L, 1)$-RW goes to a point in $D$ as the state goes to infinite then the walk is positive recurrent and the tail of the stationary distribution is geometric decay.

\subsection{{ Examples}}

Three special cases are provided here as examples.

\subsubsection{\it Degenerate to the case of $(1,1)$ state-dependent RW}

The $(1,1)$-RW is the special case in which
$q_{j}(i)\equiv 0$ for $j>1$. We denote $q_{1}(i)$ by $q(i)$. In
this case, $p(i),q(i)>0$. Let
\beqnn
    \mu_{0}=1, \quad \mu_{i}=\frac{p(0)p(1)\cdots p(i-1)}{q(1)q(2)\cdots
q(i)} \quad \mbox{for } i>0, \quad \mu=\sum_{i=0}^{\infty}\mu_{i}.
 \eeqnn
Then, Theorem~\ref{t1.1} and Theorem~\ref{t1.2} lead to the following corollaries, which are known literature results
(e.g. \cite{Chen} and \cite{Kar}).

\begin{corollary}\label{c3.1} The chain is recurrent {\it iff}~
$\sum_{i=0}^{\infty}\frac{1}{\mu_{i}p(i)}=\infty$.
\end{corollary}

\begin{corollary}\label{c3.2} The stationary distribution exists iff
$\mu<\infty$. {  In this case, $\pi_{i}=\mu_{i}/\mu$.}
\end{corollary}

For the $(1,1)$-RW, the tail of the stationary
distribution has a very simple form as follows.
\begin{corollary}\label{c3.3}
 If $p(i)\rightarrow p$ as $i\to\infty$, and $p<\frac{1}{2} $, then the stationary distribution exists and
\[
    \lim_{i\rightarrow \infty}\frac{\log\pi({i})}{i}=
-\log{\frac{1-p}{p}}.
\]
\end{corollary}

 \subsubsection{\it Degenerate to the case of $(2,1)$ state-independent RW}

We consider the classical state-independent $(2,1)$-RW on the
positive line,  for which  the transition matrix is
given by
 \beqnn
 P=\left(
             \begin{array}{cccccc}
               0 & 1 &  \\
               q_{1}+q_{2} & 0 & p & \\
               q_{2} & q_{1} & 0 & p & \\
                & q_{2} & q_{1} & 0 & p &  \\
                & & \ddots & \ddots & \ddots & \ddots \\
             \end{array}
           \right),
\eeqnn where $q_{1} +q_{2} +p =1$, and $q_{1} , q_{2} \geq 0$ and
$0<p <1$. In  this case, the recurrence criterion and
all quantities in (\ref{f1.4}), (\ref{f1.6}) and (\ref{f1.7}) can
be calculated directly from our results as follows, which lead to
a calculation of the stationary distribution.

\begin{corollary}\label{c3.4}
\beqnn
 \kappa=\infty \Longleftrightarrow q_{1}+ 2q_{2}\geq p.
\eeqnn
\end{corollary}


\begin{corollary}\label{c3.5} If $q_{1}+ 2q_{2}> p$, then the walk
$\{X_{n}\}_{n\geq 0}$ is positive recurrence, and
\beqnn
P^{i+1}[(i+1,+\infty),i-1]=\frac{\Delta-q_{1}-q_{2}}{2p}, \\
E^{i}\tau_{i}=\frac{2\Delta}{1-p-pq_{1}+3pq_{2}+(1-3p)\Delta},
\eeqnn
\beqnn
E^{i}T_{i+1}=1+\frac{1}{\lambda_1-\lambda_2}\Big(\frac{\lambda_1^{2}(1-\lambda_1^{i-1})(\lambda_2+2)}{1-\lambda_1
}-\frac{\lambda_2^{2}(1-\lambda_2^{i-1})(\lambda_1+2)}{1-\lambda_2}\Big)
+2\frac{\lambda_1^{i+1}-\lambda_2^{i+1}}{\lambda_1-\lambda_2},
\eeqnn
where $\Delta=\sqrt{(q_{1}+q_{2})^2+4pq_{2}}$, and
$\lambda_{1}$ and $\lambda_{2}$ are the eigenvalues of
\[
     M= \left(\begin{array}{cc}
\frac{q_1}{p} & \frac{q_2}{p}\\
1+\frac{q_1}{p} & \frac{q_2}{p}
\end{array} \right).
\]

\end{corollary}

\noindent {\bf Remark}~ As a consequence, we can get $\pi(i)=C(\frac{1}{\lambda_{M}})^{i}$, where $C$ is determined by the transition probabilities, and $\lambda_{M}=max\{|\lambda_{1}|,|\lambda_{2}|\}$.

\subsubsection{State-dependent queue with bulk service and
impatient customers}

This is a special case of the $(L,1)$-RW, in which $p(k)>0$ is
deceasing (arriving customers with impatience since the
probability $p(k)$ for a customer to join the system is deceasing
as the number of customers in the system increases); and the
server can simultaneously serve multiple customers up to size $L$
(bulk service). Assume $p(k) \searrow p$ and $q_j(k) \to q_j$ for
$j=1, 2, \ldots, L$. Then, according to Theorem~\ref{t1.3} if
$\sum_j jq_j > p$, then the stationary probability distribution
$\pi_k$ exists and its tail decays geometrically with rate
$1/\lambda_M$.

\vspace*{3mm}

We arrange the remainder of this paper as follows. As the main
tool, the intrinsic branching structure within the
$(L,1)$-RW will be briefly reviewed in
Section~\ref{s2} for $L=2$; and then the proofs of
the theorems and corollaries will be detailed in
Section~\ref{s4} except Theorem~\ref{t1.3}, which will be proved
in Section~\ref{s5} for a general $L$.


\section{ A brief review for the intrinsic branching structure\label{s2}}

\setcounter{equation}{0}

The intrinsic branching structure within a random walk
{ has} been studied by many authors. For the $(1,1)$-RW,
Dwass (\cite{D}, 1975) and Kesten {  \textit{et al.}
(\cite{KKS}, 1975)} observed a Galton-Watson process with the geometric offspring distribution hidden in the
nearest random walk. The branching structure is a powerful tool in
the study of random walks in a random environment
(RWRE, for short). In \cite{KKS}, Kesten \textit{et
al.}, proved a stable law for the nearest RWRE by using this
branching structure. The key point is that the hitting time
 $T_{i}$ can be calculated accurately by the branching
structure.

However, if the random walk is allowed to have jumps,
even to a bounded range, referred to as the $(L,R)$-RW, the
situation will become much more complicated.

For the $(L,1)$-RW, when the walk starts at $0$ and
$\limsup_{n\rightarrow \infty}X_n= \infty$, a multi-type branching
process has been revealed by Hong \textit{et al.}(\cite{hw}, 2009). for calculating the hitting time $T_{1}=\inf\{n>0,
X_{n}>0\}$. A similar work has been done for the $(1,R)$-RW
(\cite{hzh}, 2010). It must be emphasized that these two branching structures are not symmetric, instead
they are essentially different.

{For the purpose of calculating the stationary
distribution in this paper, both structures will be used. More
specifically,  for calculating $E^{i}T_{i}$ for $i>0$, if the
first step is down (possible at $i-1$ or $i-2$ when $L=2$), the
branching structure within the $(L,1)$-RW is used,  and
alternatively if the first step is up (to state $i+1$), the
branching structure within the $(1,R)$-RW is used. Here we call
them the ``lower" and the ``upper" branching structures
respectively, which we will briefly introduce below. Note that if
we assume $q_{2}(i)\equiv 0$, both branching structures degenerate
to the case of the $(1,1)$-RW.

\subsection{The ``lower" branching structure}
The following discussion is based on $L=2$. The
general case can be similarly discussed, which is much more
complicated. Assume that $X_{0}=i$, if the first step is down
(possible at $i-1$ or $i-2$) we can calculate $E^{i}T_{i}$ by
using the branching structure within the $(L,1)$-RW (\cite{hw},
2009). Define
\[
    U^{1}_{0}=\#\{0<j<T_{i+1}:X_{j-1}=1,X_j=0\} \quad \text{ and } \quad U^{2}_{0}=0;
\]
and
\[
    U^{l}_{k}=\#\{0<j<T_{i+1}:X_{j-1}>k,X_j=k-l+1\} \quad \text{
    for } 1\leq k<i+1, \; l=1,2.
\]
Setting
\[
    U_k=(U^{1}_{k},U^{2}_{k}) \quad \text{ for } 0\leq k<i+1.
\]
We then have the following property.
\medskip

\noindent{\bf Theorem A} (Hong and Wang~\cite{hw})
{\it\noindent(1) The process $\{U_k\}_{k=i}^{0}$}
is a $2$-type branching process whose branching mechanism is given
by:
\beqnn
&&P(U_{0}=(a,0)\big|U_{1}=e_1)=\big(q_{1}(1)+q_{2}(1)\big)^{a}p(1),\\
&&P(U_{0}=(1+a,0)\big|U_{1}=e_2)=\big(q_{1}(1)+q_{2}(1)^{a}p(1);
\eeqnn
 and for $k>1$,
  \beqnn
&&P(U_{k-1}=(a,b)\big|U_{k}=e_1)=\frac{(a+b)!}{a!
b!}q_{1}(k)^{a}q_{2}(k)^{b}p(k),\\
&&P(U_{k-1}=(1+a,b)\big|U_{k}=e_2)=\frac{(a+b)!}{a!
b!}q_{1}(k)^{a}q_{2}(k)^{b}p(k).
\eeqnn

{\it (2) For the process $\{U_k\}_{k=i}^{0},$ let $ M_k$ be
the $2\times 2$ mean matrix whose $l$-th row is
$E(U_{k-1}|U_k=e_l).$ Then, one has that
 \beqnn M_{1}= \left(
\begin{array}{cc}
\frac{q_{1}(1)+q_{2}(1)}{p(1)} & 0\\
\frac{1}{p(1)} & 0
\end{array} \right), \quad
M_{k}= \left(
\begin{array}{cc}
\frac{q_{1}(k)}{p(k)} & \frac{q_{2}(k)}{p(k)}\\
1+\frac{q_{1}(k)}{p(k)} & \frac{q_{2}(k)}{p(k)}
\end{array} \right), \quad k>1.
 \eeqnn

(3) { $E^{0}T_{1}=1$, and for $i >0$,
$T_{i+1}=1+\sum_{k=0}^{i}U_k\cdot(2,1)'$
 and}
  \beqlb\label{f1.7}
 E^{i}T_{i+1}=1+\sum_{k=0}^{i}EU_k\cdot(2,1)'=1+\sum_{k=1}^{i}e_{1} M_{i}M_{i-1}\cdots M_{i-k+1}(2,1)'.
 \eeqlb}

 \noindent{\bf Remark}  i) The positions of the walk correspond
to the time of the branching process. For example, { in
our notation,} $U_{k}$ is indeed the $(i-k)$-th generation of the
branching process.

ii) The condition $\limsup_{n\rightarrow \infty}X_n=+\infty$ in
\cite{hw} is obviously satisfied in our reflecting model.

iii) It is not difficult to understand that
{ the branching structure from} $U_{1}$ to $U_{0}$ is
different from others because of reflecting. We omit the proof
here.

\subsection{The ``upper" branching structure}

Assume that $X_0=i$, $i>0$. If the first step is up (at $i+1$), we
can calculate $E^{i+1}T_{i}$ by using the branching structure
within $(1,R)$-RW (\cite{hzh}, 2010). If the $(2,1)$-reflecting
random walk is recurrent, $\tau_{k}$, $k\geq i$, is defined
 the same as in Section~$1$. Note that $\tau_k<\infty$
$P$-a.s. To calculate $\tau_i$ accurately, Hong and Zhang
(\cite{hzh}, 2010) defined a multi-type branching
process by decomposing the path of the walk. Intuitively,  if the
walk from $k\geq i$ takes a step to $k+1$, it must
 across back to $k$ or jump over $k$ (to $k-1$)
because of $\tau_k<\infty$ $P$-a.s., in which there are only three
  ways of moving down: from $k+1$ to $k$, from $k+2$
to $k$ and from $k+1$ to $k-1$. So we divide all the steps from
$k$ to $k+1$ into three  kinds of steps according to
the above three ways of moving down. Let $A(k)$, $B(k)$ and
$C(k)$ be the numbers of steps from $k$ to $k+1$ before time
$\tau_i$  corresponding to moving from $k+1$ to $k$,
from $k+2$ to $k$ and from $k+1$ to $k-1$, respectively.
 As for the last step of $\tau_i$, we can consider it
as a immigration for the multi-type branching processes.

Consider integers $n\geq i> 1$, and define the exit probabilities:
\beqnn
P^i[(i,n),i-1]= P^i \{\mbox{$X_{n}$ leaving $\{i,i+1,\ldots,n-1,n\}$ at the point $i-1$}\},\\
P^i[(i,n),i-2]= P^i \{\mbox{$X_{n}$ leaving
$\{i,i+1,\ldots,n-1,n\}$ at the point $i-2$}\}.
 \eeqnn
 In  Hong and Zhang~\cite{hzh} { (see Lemma 2.1)}, it has been calculated that
(see also \cite{B02})
\beqlb\label{f1.4}
    P^{i}[(i,n),i-1]=\frac{\langle e_1,[\tilde M_{i}+\cdots+\tilde M_{n}\cdots \tilde
M_{i}]v\rangle}
{1+\langle e_1, [\tilde M_{i}+\cdots+\tilde M_{n}\cdots \tilde M_{i}]e_1\rangle},\nonumber\\
P^{i}[(i,n),i-2]= \frac{\langle e_1,[\tilde M_{i}+\cdots+\tilde
M_{n}\cdots \tilde M_{i}]e_2\rangle} {1+\langle e_1,[\tilde
M_{i}+\cdots+\tilde M_{n}\cdots \tilde M_{i}]e_1\rangle},
 \eeqlb where $v=e_1'-e_2'$, and $\tilde M_{i}= \left(
\begin{array}{cc}
\frac{q_1(i)+q_2(i)}{p(i)} & \frac{q_2(i)}{p(i)}\\
1 & 0
\end{array} \right),~i\geq 1.$\\

\noindent If $\kappa=\infty$, let \beqlb\label{f1.5}
& &\gamma(i)=p(i)\cdot P^{i+1}[(i+1,+\infty),i-1],\nonumber\\
& &\alpha(i)=p(i)\cdot P^{i+1}[(i+1,+\infty),i]\cdot\frac{q_1(i+1)}{q_1(i+1)+\gamma(i+1)},\nonumber\\
& &\beta(i)=p(i)\cdot
P^{i+1}[(i+1,+\infty),i]\cdot\frac{\gamma(i+1)}{q_1(i+1)+\gamma(i+1)}.
\eeqlb

Set for $k\geq i$,
\[
     V(k)=[A(k),B(k),C(k)].
\]
 Then we have the following branching structure
within the $(2,1)$-RW.
\medskip

\noindent{\bf Theorem B} ({ Hong and Zhang~\cite{hzh}})
{\it Assume $\kappa=\infty$. Then,

(1) $\Big(V(k)=[A(k),B(k),C(k)]\Big)_{k\geq i}$ is an
inhomogeneous multi-type branching process with immigration \beqnn
V(i-1)=[1,0,0],~~\mbox{with probability}~~\frac{q_1(i)}{1-\alpha(i)-\beta(i)},\\
V(i-1)=[0,1,0],~~\mbox{with probability}~~\frac{\gamma(i)}{1-\alpha(i)-\beta(i)},\\
V(i-1)=[0,0,1],~~\mbox{with probability}~~\frac{q_2(k)}{1-\alpha(i)-\beta(i)}.
\eeqnn
The offspring distribution is given by
\beqnn
& &P\Big(V(k)=[a,b,0]~\Big|~V(k+1)=[1,0,0]\Big)=[1-\alpha(k)-\beta(k)]C_{a+b}^a\alpha(k)^a\beta(k)^b,\\
& &P\Big(V(k)=[a,b,1]~\Big|~V(k+1)=[0,1,0]\Big)=[1-\alpha(k)-\beta(k)]C_{a+b}^a\alpha(k)^a\beta(k)^b,\\
&
&P\Big(V(k)=[a,b,0]~\Big|~V(k+1)=[0,0,1]\Big)=[1-\alpha(k)-\beta(k)]C_{a+b}^a\alpha(k)^a\beta(k)^b,
\eeqnn where $\alpha(k)$, $\beta(k)$ and $\gamma(i)$ are defined
in (\ref{f1.4}).

(2) The offspring mean matrix of the { $(k-i+1)$-st}
generation of the multi-type branching process is \beqnn N_{k}=
\left(
\begin{array}{ccc}
\frac{\alpha(k)}{1-\alpha(k)-\beta(k)} & \frac{\beta(k)}{1-\alpha(k)-\beta(k)} &0\\
\frac{\alpha(k)}{1-\alpha(k)-\beta(k)} & \frac{\beta(k)}{1-\alpha(k)-\beta(k)} &1\\
\frac{\alpha(k)}{1-\alpha(k)-\beta(k)} & \frac{\beta(k)}{1-\alpha(k)-\beta(k)} &0
\end{array} \right).
\eeqnn

(3)
$\tau_{i}=1+\sum_{k=i}^{\infty}\big[2A(i)+2B(i)+C(i)\big]=1+\big<(2,2,1),\sum_{k=i}^{\infty}V(i)\big>$,
and \beqlb\label{f1.6}
E^{i}\tau_{i}=1+\Big\langle(2,2,1),~\frac{1}{1-\alpha(i)-\beta(i)}\Big(q_1(i),~\gamma(i),~q_2(i)\Big)
\cdot \sum_{k=1}^{+\infty}N_{1}\cdots N_{k}\Big\rangle. \eeqlb }

\noindent{\bf Remark} i) In \cite{hzh}, Hong and Zhang considered
the branching structure within the $(1,2)$-RW on the line starting at $0$ before the ladder time $\inf\{k>0,
X_{k}>0\}$. It corresponds to the ``upper" part of
our model.

ii) $E^{1}\tau_{1}<\infty $ is the sufficient condition of
$\kappa=\infty$.






 \section{Proofs\label{s4}}

 \setcounter{equation}{0}

  \noindent{\it Proof of Theorem 1.1} ~~
According to Fact~1, we can prove the result by two
different methods. Here we prove it through
calculating $E^0N(0)$ directly by using the branching
structure; and in the appendix we provide an
analytical proof by solving the system of infinite linear
equations. Actually, we find the solution to the system of
infinite linear equations also in terms of the observation of the
branching structure.

Recall $N(0)=\sum_{n=1}^{\infty} 1_{\{X_n=0\}}$ is the occupation
time of position 0. Because the walk is reflected at 0 with
probability 1, the walk goes to $+\infty$ a.s.  Note
that the walk goes up skip freely. Therefore, we can decompose
the whole trajectory of the walk as the combination of the pieces
from position $i$ to $i+1$ for $i=0, 1, 2, \ldots$,}
i.e., $\{X_n, n\geq 0\}=\bigcup_{i=0}^{\infty} \{X_n, T_i\leq n<
T_{i+1}\}$, where $T_0=0$. So, to calculate $N(0)$
the occupation time of position 0, we need only to calculate the
occupation time of position 0 in each piece $\{X_n,
T_i\leq n< T_{i+1}\}$ of the trajectory, in which  a ``lower"
multi-type process with an immigration at position
$i$ is hidden (Theorem A). To this end,
  let  $\theta_{1}=1$ and for $ i>1$,
\beqnn
& &\theta_{i}=\sharp \{k: T_{i-1}\leq k< T_{i}, X_{k}=0\};
\eeqnn

\

\

 \begin{center}
  \includegraphics[scale=.4,angle=-90,totalheight=60mm]{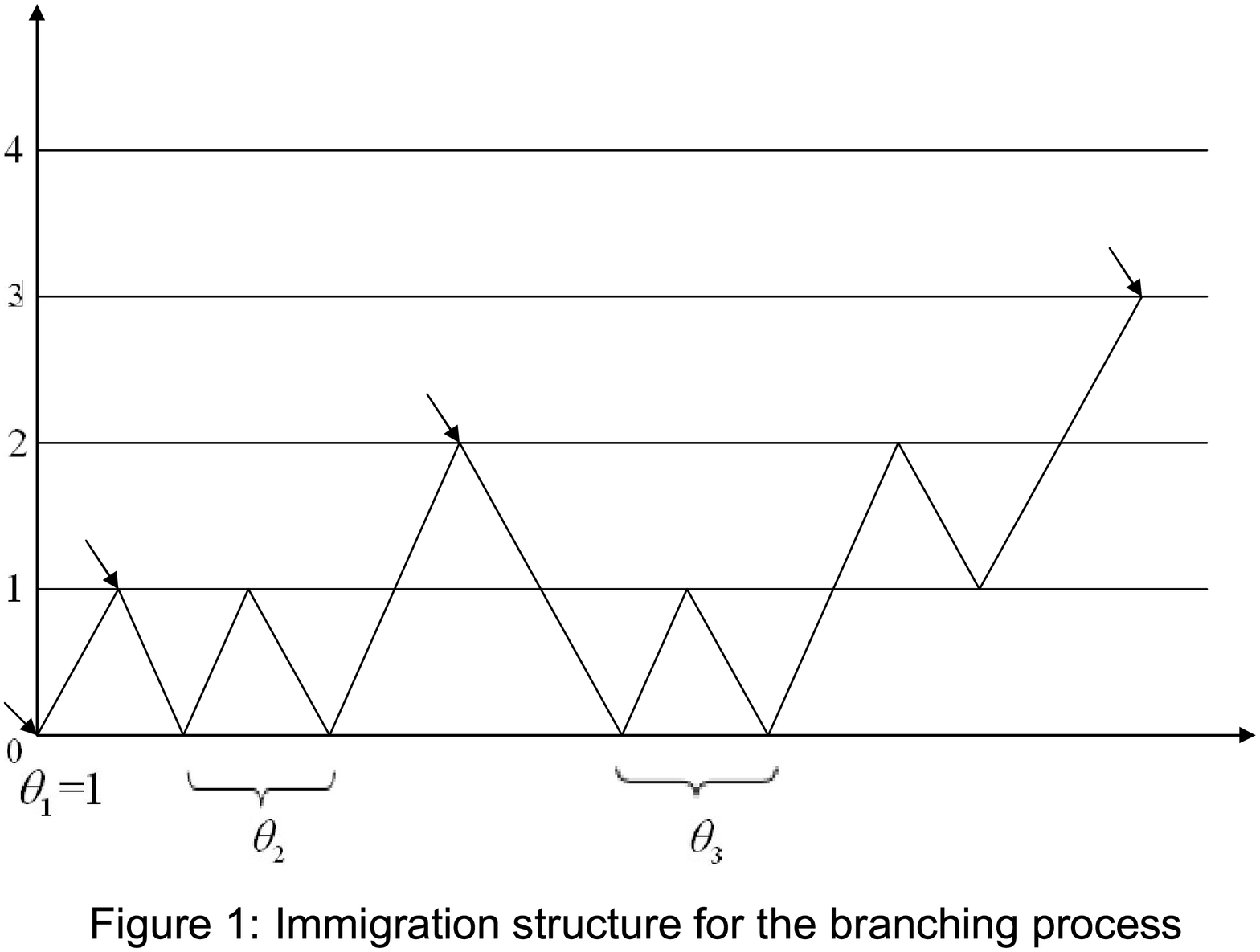}
\end{center}
which is the population of the $i$-th generation of the  multi-type branching process with a immigration at $i$. By Theorem A, we have
\beqnn
& &E^{0}\theta_{i}=e_{1}M_{i-1}M_{i-2}\cdots M_{1}u,~~~i\geq 1.
 \eeqnn

From the decomposition we know that
$N(0):=\sum_{n=1}^{\infty} 1_{\{X_n=0\}}=\sum_{i=1}^{\infty}\theta_{i}.$ As a consequence we get
$$
E^{0}N(0)=E^{0} \sum_{i=1}^{\infty}\theta_{i}=\sum_{i=1}^{\infty}e_{1} M_{i-1} M_{i-2}\cdots M_{1}u:=\kappa.
$$
The proof is complete by Fact~1. \qed

\noindent{\bf Remark} Let $\xi_{0}=0,$ $\xi_{1}=1$ and
$\xi_{n}=\sum_{i=1}^{n}\theta_{i}$ for  $n\geq 1$,  the
occupation time at position 0 before the walk hitting position
$n$. If $y_n=E^{0}\xi_{n}$, we find that $y_n$ is a
solution of the system of infinite linear equations
$\sum_{j=0}^{+\infty}P_{ij}y_{j}=y_{i} ~i\geq0$. See the appendix,
in which an analytical version proof of the present
theorem is provided.

\medskip
\noindent{\it Proof of Theorem 1.2}~~In our model, $X_{n}$ is
irreducible  and note that $E^0T_0=1+E^{1}\tau_{1}$, where the
explicit expression of $E^i\tau_i$ is given in (\ref{f1.6}). So
$E^{1}\tau_{1}<\infty$ assures that the walk $X_{n}$
is positive recurrent by Fact~2 and has the stationary
distribution $\pi$ with $\pi(i)=\frac{1}{E^{i}T_{i}}$. To
calculate $E^{i}T_{i}$, we consider the first step of
the walk starting at $X_0=i$. There are four possible types of
the trajectory of the walk from position $i$ back to position $i$
as shown in the following graph.

\

\

\begin{center}
  \includegraphics[scale=.4,angle=-90,totalheight=60mm]{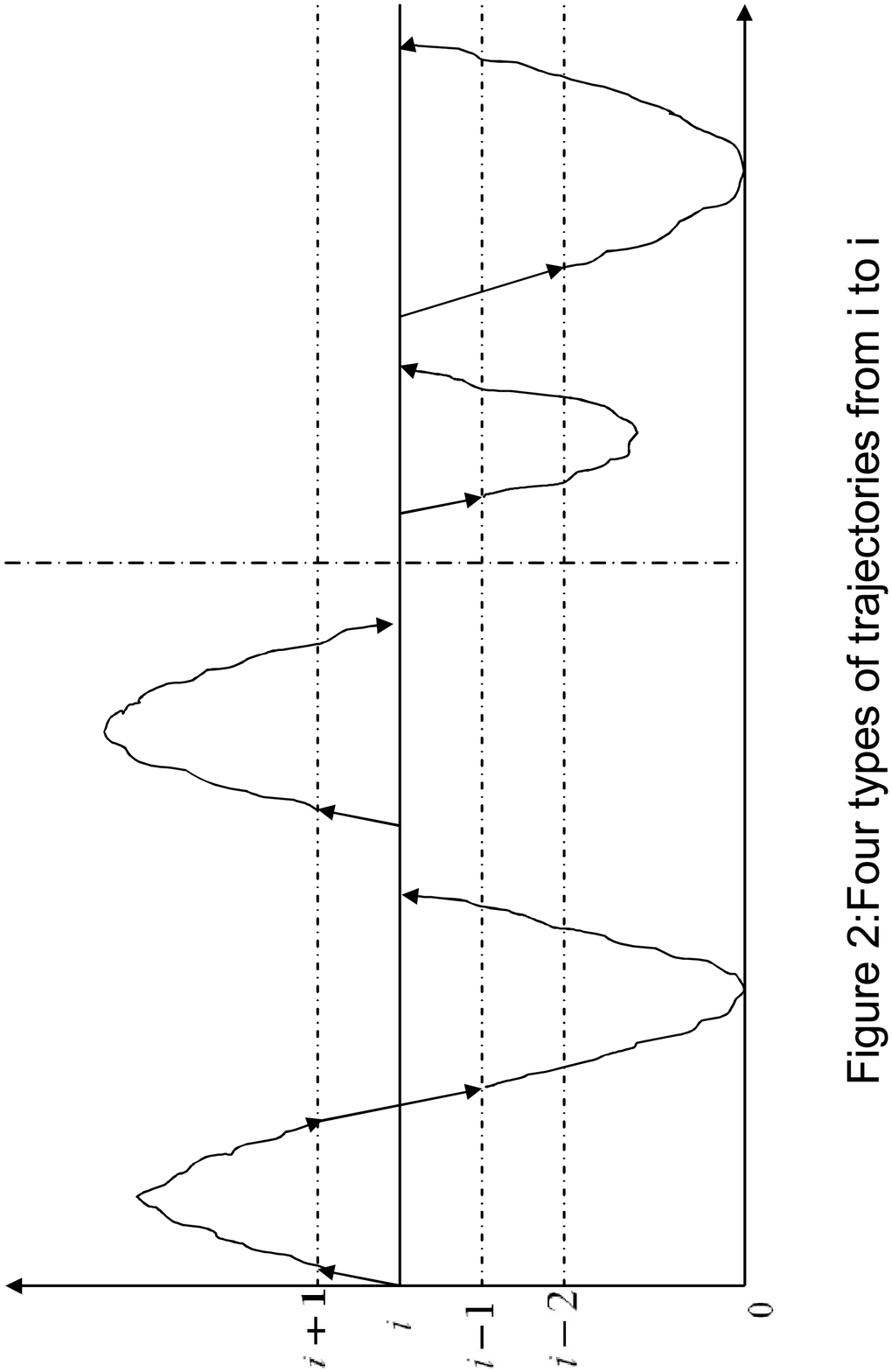}
\end{center}

By the Markov property  one can get, for $i>0$,
\beqlb\label{ff}
E^{i}T_{i}&=&p(i)(E^{i+1}T_{i}+1)+q_{1}(i)(E^{i-1}T_{i}+1)+q_{2}(i)(E^{i-2}T_{i}+1) \nonumber\\
&=&p(i)(E^{i+1}T_{i}+1)+q_{1}(i)(E^{i-1}T_{i}+1)+q_{2}(i)(E^{i-2}T_{i-1}+E^{i-1}T_{i}+1)\nonumber\\
&=&p(i)E^{i+1}T_{i}+(q_{1}(i)+q_{2}(i))E^{i-1}T_{i}+q_{2}(i)E^{i-2}T_{i-1}+1,
\eeqlb and $E^{0}T_{0}=E^{1}T_{0}+1$.

We now calculate $E^{i+1}T_{i}$. For the walk
starting at $i+1$, there are only two kinds of ways to hit $i$
for the first time: from $i+1$ or $i+2$ to $i$ with probability
$P^{i+1}[(i+1,+\infty),i]$; or first from $i+1$  to $i-1$, then
from $i-1$ to $i$ with the probability
$P^{i+1}[(i+1,+\infty),i-1]$. Consequently,
 \beqnn
    E^{i+1}T_{i}&=& P^{i+1}[(i+1,+\infty),i]E^{i+1}\tau_{i+1}+P^{i+1}[(i+1,+\infty),i-1](E^{i+1}\tau_{i+1}+E^{i-1}T_{i})\\
&=&E^{i+1}\tau_{i+1}+P^{i+1}[(i+1,+\infty),i-1]E^{i-1}T_{i}.
 \eeqnn
By $(\ref{ff})$, we get
 \beqnn
 E^{i}T_{i}
 =p_{i}E^{i+1}\tau_{i+1}+(q_{1}(i)+q_{2}(i)+p(i)P^{i+1}[(i+1,+\infty),i-1)]E^{i-1}T_{i}+q_{2}(i)E^{i-2}T_{i-1}+1.
 \eeqnn \qed

\medskip
\noindent{\it  Proof of Corollary \ref{c3.1} } ~~
In this case, we have for $i\geq 1$, $M_{i}:= \left(
\begin{array}{cc}
\frac{q(i)}{p(i)} & 0\\
\frac{1}{p(i)} & 0
\end{array} \right).$
Therefore,
\beqnn
\kappa=\sum_{k=1}^{\infty}e_{1}M_{k}M_{k-1}\ldots M_{1}u=\sum_{k=1}^{\infty}(1,~0)\left(
\begin{array}{cc}
\frac{q_{k}q_{k-1}\ldots q_{i}}{p_{k}p_{k-1}\ldots p_{i}} & 0\\
\frac{q_{k-1}q_{k-2}\ldots q_{i}}{p_{k}p_{k-1}\ldots p_{i}} & 0
\end{array} \right)(1,~1)'=\sum_{i=1}^{\infty}\frac{1}{\mu_{i}p_{i}}.
\eeqnn
We also have $\frac{1}{\mu_{0}p_{0}}=1.$ As a result, $\kappa=\infty \Leftrightarrow \sum_{i=0}^{\infty}\frac{1}{\mu_{i}p_{i}}=\infty.$
\qed

\

\noindent{\it  Proof of Corollary \ref{c3.2} }~~ By (\ref{f1.4}), (\ref{f1.5}), (\ref{f1.6}) and (\ref{f1.7}), one can get
$$E^{i}\tau_{i}=\frac{1}{q(i)}+\frac{p(i)}{q(i)q(i+1)}+\frac{p(i)p(i+1)}{q(i)q(i+1)q(i+2)}+\cdots~,~~~~i\geq 1,$$
$$E^{i}T_{i+1}=\frac{1}{p(i)}+\frac{q(i)}{p(i-1)p(i)}+\frac{q(i-1)q(i)}{p(i-2)p(i-1)p(i)}+\cdots+
\frac{q(1)q(2)\ldots q(i)}{p(0)p(1)\ldots p(i)},~~~~i\geq 0,$$
So $\mu=E^{1}\tau_{1}$, and $\mu<\infty \Longleftrightarrow E^{1}\tau_{1}<\infty$. By Theorem 1.2,
\beqnn
\pi(i)=\frac{1}{p(i)E^{i+1}\tau_{i+1}+q(i)E^{i-1}T_{i}+1}=\frac{\mu_{i}}{\mu}.
\eeqnn
\qed

\medskip
\noindent{\it  Proof of Corollary \ref{c3.4} }~~
We define
\beqnn M_{1}= \left(
\begin{array}{cc}
\frac{q_{1}+q_{2}}{p} & 0\\
\frac{1}{p} & 0
\end{array} \right):=A,~~~
M_{i}= \left(
\begin{array}{cc}
\frac{q_{1}}{p} & \frac{q_{2}}{p}\\
1+\frac{q_{1}}{p} & \frac{q_{2}}{p}
\end{array} \right):=B,~~~i>1
\eeqnn
\beqnn
\left(
\begin{array}{cc}
1 & 0 \\
1 & 1 \\
\end{array}\right)\cdot
\left(
\begin{array}{cc}
\lambda_1+\lambda_2 & 0\\
1 & 1 \\
\end{array}
\right):=C,~~~
\left(
\begin{array}{cc}
1 & 0 \\
1 & 1 \\
\end{array}\right)\cdot
\left(
\begin{array}{cc}
\lambda_1 & \lambda_2\\
1 & 1 \\
\end{array}
\right):=D.
\eeqnn

Let $\lambda_{1,2}=\frac{q_1+q_2\pm\sqrt{(q_1+q_2)^2+4pq_2}}{2p}$
be the eigenvalues of $M_{i},i>1$,and let $\tilde \lambda_{1,2}$ be the eigenvalues of $M_{1}$. We can see that $\tilde \lambda_{1}=\lambda_{1}+\lambda_{2}$ and $\tilde \lambda_{2}=0$. Then, we can decompose $M_{i}$ as
$$
M_{1}=C\cdot
\left(
\begin{array}{cc}
\lambda_1 & 0 \\
0 & \lambda_2 \\
\end{array}\right)
\cdot C^{-1},~~~
M_{i}=D\cdot
\left(
\begin{array}{cc}
\lambda_1 & 0 \\
0 & \lambda_2 \\
\end{array}\right)
\cdot D^{-1},~~~i>1
$$
Then, if $\lambda_1\neq 1$
\beqnn
\nonumber & &\sum_{k=1}^{n}e_{1}M_{k} M_{k-1}\ldots M_{1}u=e_{1}(E+B+B^{2}\ldots B^{n})Au\\
&=&
e_{1}\cdot D
\cdot
\left(
\begin{array}{cc}
\frac{1-\lambda_1^{n}}{1-\lambda_1} & 0 \\
0 & \frac{1-\lambda_2^{n}}{1-\lambda_2} \\
\end{array}\right)
\cdot
D^{-1}\cdot C
\cdot
\left(
\begin{array}{cc}
\lambda_1+\lambda_2 & 0 \\
0 & 0 \\
\end{array}\right)
\cdot
C^{-1}\cdot(1,1)'\\
&=&\frac{1}{\lambda_1-\lambda_2}\Big(\frac{1-\lambda_1^{n}}{1-\lambda_1}\lambda_1^{2}-\frac{1-\lambda_2^{n}}{1-\lambda_2}\lambda_2^{2}\Big),
\eeqnn
if $\lambda_1=1$
\beqnn
\nonumber \sum_{k=1}^{n}e_{1}M_{k}M_{k-1}\ldots M_{1}u =\frac{1}{\lambda_1-\lambda_2}\Big(\lambda_1^{2}n-\frac{1-\lambda_2^{n}}{1-\lambda_2}\lambda_2^{2}\Big). \\
\eeqnn
It is easy to see that $\lambda_{2}=\frac{q_1+q_2-\sqrt{(q_1+q_2)^2+4pq_2}}{2p}\in (-1,0]$. We get $\lambda_{2}^{n}\rightarrow 0$ as $n\rightarrow \infty$.

Hence, $\sum_{k=1}^{\infty }e_{1}M_{k}M_{k-1}\ldots M_{1}u=\infty \Leftrightarrow \lambda_{1}\geq 1.$ By some calculation,
\beqnn
\lambda_{1}\geq 1\Leftrightarrow q_{1}+ 2q_{2}\geq p.
\eeqnn \qed

\medskip
\noindent{\it  Proof of Corollary \ref{c3.5} } ~~
Define $A,B,C,D,\lambda_{1}$ and $\lambda_{2}$ the same as before. $P^{i+1}[(i+1,+\infty),i-1]$ and $E^{i}\tau_{i}$ are calculated in \cite{hzh}, we need only to calculate $E^{i}T_{i+1}$.
By   Theorem A, we   get
\beqnn
E^{i}T_{i+1}
&=&1+e_{1}M_{i}(2,1)'+e_{1}M_{i} M_{i-1}(2,1)'+\ldots+e_{1}M_{i}M_{i-1}\ldots M_{1}(2,1)'\\
&=&1+e_{1}(B+B^{2}+\ldots+B^{i-1})(2,1)'+e_{1}(B_{i-1}A)(2,1)'\\
&=&1+e_{1}\cdot D
\cdot
\left(
\begin{array}{cc}
\frac{\lambda_1(1-\lambda_1^{i-1})}{1-\lambda_1} & 0 \\
0 & \frac{\lambda_2(1-\lambda_2^{i-1})}{1-\lambda_2} \\
\end{array}\right)
\cdot
D^{-1}\cdot(2,1)'\\
&
&+
e_{1}\cdot D
\cdot
\left(
\begin{array}{cc}
\lambda_1^{i-1} & 0 \\
0 & \lambda_2^{i-1} \\
\end{array}\right)
\cdot
D^{-1}
\cdot C\cdot
\left(
\begin{array}{cc}
\lambda_1+\lambda_2&0 \\
1 & 1 \\
\end{array}
\right)
\cdot
\left(
\begin{array}{cc}
\lambda_1+\lambda_2 & 0 \\
0 & 0 \\
\end{array}\right)
\cdot
C^{-1}\cdot(2,1)'\\
&=&1+\frac{1}{\lambda_1-\lambda_2}(\frac{\lambda_1^{2}(1-\lambda_1^{i-1})(\lambda_2+2)}{1-\lambda_1
}-\frac{\lambda_2^{2}(1-\lambda_2^{i-1})(\lambda_1+2)}{1-\lambda_2})
+2\frac{\lambda_1^{i+1}-\lambda_2^{i+1}}{\lambda_1-\lambda_2}.
\eeqnn \qed

\section{Tail asymptotic of $\pi (i)$ --- proof of Theorem \ref{t1.3} \label{s5}}

\setcounter{equation}{0}

It is well  known that when the random walk is
state-independent the tail of the stationary distribution of the
walk has a geometric decay. With the help of the
explicit expression of the stationary distribution, we can
consider the tail behavior in  the state-dependent case. We note
that the expression of the stationary distribution in
Theorem~\ref{t1.2} is given in terms of the
decomposition of the trajectory with different
parts. This enables us to find out that the key factor to
determine the tail asymptotic behavior of $\pi (i)$ is
$E^{i}T_{i+1}$ ( because of $E^{1}\tau_{1}<\infty $). In this
section, we consider the case with a general $L$ (not
only for $L=2$), i.e., the $(L,1)$-RW.

Suppose $L\geq1$, to express $E^{i}T_{i+1}$ we need only the ``lower" branching structure  hidden in the random walk which we have introduced in Section 2.  Recall for $i>0$, $$
 E^{i}T_{i+1}=1+\sum_{k=0}^{i}EU_k\cdot w_{L}'=1+\sum_{k=1}^{i}e_{1} M_{i}M_{i-1}\ldots M_{i-k+1}w_{L}',
$$where $e_{L}=(1,0,\cdots,0)$, $w_{L}=(2,1,\cdots,1)$. The
offspring mean matrices of the multi-type branching
processes are given by:
 \beqnn M_{i}=\left(\begin{array}{ccccc}
\frac{q_{1}(i)}{p(i)} & \frac{q_{2}(i)}{p(i)} & \cdots & \frac{q_{L-1}(i)}{p(i)} & \frac{q_{L}(i)}{p(i)}\\
1+\frac{q_{1}(i)}{p(i)} & \frac{q_{2}(i)}{p(i)} & \cdots & \frac{q_{L-1}(i)}{p(i)} & \frac{q_{L}(i)}{p(i)}\\
\vdots & \vdots & \ddots & \vdots & \vdots\\
\frac{q_{1}(i)}{p(i)} & \frac{q_{2}(i)}{p(i)} &\cdots & \frac{q_{L-1}(i)}{p(i)} & \frac{q_{L}(i)}{p(i)}\\
\frac{q_{1}(i)}{p(i)} & \frac{q_{2}(i)}{p(i)} & \cdots & 1+\frac{q_{L-1}(i)}{p(i)} & \frac{q_{L}(i)}{p(i)}
\end{array} \right)_{L\times L},~i\geq L,~~
\eeqnn
with a little attention to the reflect effect we notice the difference for $1\leq i<L$,
\beqnn
M_{L-1}=\left(\begin{array}{ccccc}
\frac{q_{1}(i)}{p(i)} & \frac{q_{2}(i)}{p(i)} & \cdots & \frac{q_{L-1}(i)+q_{L}(i)}{p(i)} & 0\\
1+\frac{q_{1}(i)}{p(i)} & \frac{q_{2}(i)}{p(i)} & \cdots & \frac{q_{L-1}(i)+q_{L}(i)}{p(i)} & 0\\
\vdots & \vdots & \ddots & \vdots & \vdots\\
\frac{q_{1}(i)}{p(i)} & \frac{q_{2}(i)}{p(i)} &\cdots & \frac{q_{L-1}(i)+q_{L}(i) }{p(i)} & 0\\
\frac{q_{1}(i)}{p(i)} & \frac{q_{2}(i)}{p(i)} & \cdots & 1+\frac{q_{L-1}(i)+ q_{L}(i)}{p(i)} & 0
\end{array} \right)_{L\times L}, ~~ \ldots \ldots,
\eeqnn
and
\beqnn
M_{1}=\left(\begin{array}{ccccc}
\frac{q_{1}(i)+\cdots +q_{L}(i)}{p(i)} & 0 & \cdots & 0 & 0\\
1+\frac{q_{1}(i)+\cdots +q_{L}(i)}{p(i)} & 0 & \cdots & 0 & 0\\
\vdots & \vdots & \ddots & \vdots & \vdots\\
\frac{q_{1}(i)+\cdots +q_{L}(i)}{p(i)} &  0 & \cdots & 0 & 0\\
\frac{q_{1}(i)+\cdots +q_{L}(i)}{p(i)} &  0 & \cdots & 0 & 0
\end{array} \right)_{L\times L}.
\eeqnn

The following two propositions about matrix analysis
are needed for our proof. Recall, for a matrix $A$, $\rho(A)$ is
the spectral radius of $A$, and $\lambda_{A}$ is the maximum
eigenvalue of it.

\begin{proposition} (Perron's Theorem, \cite{hor})
If $A=(a_{ij})_{L\times L}\in \mathbb{R}^{L\times L}$. and $A>0$(which means $\forall i,j>0, a_{ij}>0$), then\\
(a)$\rho(A) >0$;\\
(b) $\rho(A)$ is an eigenvalue of A, and it is the unique eigenvalue of maximum modulus;\\
(c) $\rho(A)$ is an algebraically (and hence geometrically) simple
value of A;\\
(d) $[\rho(A)^{-1}A]^{k}\rightarrow R$ as $k\rightarrow \infty$, where $R\in \mathbb{R}^{L\times L}$, and $R>0$. \qed
\end{proposition}

A measure for the distance between the spectra
$\sigma(A)=\{\lambda_{1},\lambda_{2},\ldots \lambda_{L}\}$
and $\sigma(B)=\{\mu_{1},\mu_{2},\ldots \mu_{L}\}$ is
defined below, which is the optimal matching
distance:
$$d(\sigma(A),\sigma(B))=\min\limits_{\theta\in S_{L}}\max\limits_{i\in\{1,2,\cdots L\}}|\lambda_{i}-\mu_{\theta_{i}}|,$$
where $S_{L}$ denotes the group of all permutations of the numbers $1,2,\cdots L$

\begin{proposition} (\cite{Kra},\cite{Ost})  
$$d(\sigma(A),\sigma(B))\leq 4(2K)^{1-\frac{1}{n}}\|A-B\|^{\frac{1}{n}}$$
where $K=\max\{\|A\|,\|B\|\}$.
\end{proposition}

This result says that there exists a permutation such
that the maximum of the distance between the corresponding
eigenvalues is small enough.

\

Let $D:=\{(p,q_{1},q_{2},\cdots q_{L}): p+\sum_{j=1}^{L}q_{j}=1; ~~ \sum_{j=1}^{L}jq_{j}>p;~~\forall j, ~q_{j}\geq 0, p>0 \}.$ A point $P=(p,q_{1},q_{2},\cdots q_{L})$ in $D$ corresponds a state-independent transition probability of the $(L,1)$-RW, and the offspring mean matrix of the correspond multi-type process is denoted as $M$ (it is independent of the position $i$).

\begin{lemma}
If $P\in D$, the offspring mean matrix of the lower branching process $M=(m_{ij})_{L\times L}$ has the unique eigenvalue of maximum modulus. That is
$$\rho(M)>1,$$
and $\rho(M)=\lambda_{M}$. Let $C=4(\frac{2}{p}+2L)^{1-\frac{1}{L}}L^{\frac{1}{L}}$.  For $\varepsilon<\min\{\frac{1}{C^{L}}(\lambda_{M}-1)^{L},\min \limits_{i,j}m_{ij}\}$,
we have
$$\lambda_{M}^{\varepsilon-}:=\rho(M-\varepsilon E)>\lambda_{M}-C\varepsilon^{\frac{1}{L}}>1 , ~~~~ 1<\lambda_{M}^{\varepsilon+}:=\rho(M+\varepsilon E)<\lambda_{M}+C\varepsilon^{\frac{1}{L}} ,$$
\end{lemma}
where $E=(1,1,\cdots 1)'\cdot(1,1,\cdots 1).$

\medskip
\noindent{\it Proof }~~
First, to prove $\rho(M)>1$,
let $  \hat{M}=\left(\begin{array}{cc}
1 & 0\\
0 & M
\end{array} \right)$,~
by some calculations, we get
$$|\lambda I-\hat{M}|=
\left(\begin{array}{ccccc}
\lambda-1 & 0 & \cdots & 0 & 0\\
0 & \lambda-\frac{q_{1}}{p} & \ldots & -\frac{q_{L-1}}{p} & -\frac{q_{L}}{p}\\
\vdots & \vdots & \ddots & \vdots & \vdots\\
0 & -1-\frac{q_{1}}{p} & \cdots & \lambda-\frac{q_{L-1}}{p} & -\frac{q_{L}}{p}\\
0 & -\frac{q_{1}}{p} & \cdots & -1-\frac{q_{L-1}}{p} & \lambda-\frac{q_{L}}{p}
\end{array} \right)
=
\left(\begin{array}{ccccc}
\lambda-\frac{1}{p} & \frac{q_{1}}{p} & \cdots & \frac{q_{L-1}}{p} & \frac{q_{L}}{p}\\
-1 & \lambda & \ldots & 0 & 0\\
\vdots & \vdots & \ddots & \vdots & \vdots\\
0 & 0 & \cdots & \lambda & 0\\
0 & 0 & \cdots & -1 & \lambda
\end{array} \right)$$.

Define $f(x)=x^{L+1}-\frac{1}{p}x^{L-1}+\cdots+\frac{q_{L-1}}{p}x+\frac{q_{L}}{p}$.
We have $|\lambda I-\hat{M}|=f(x).$ Let $F(x)=f(\frac{1}{x})$, then,
$$F(x)=q_{L}x^{L+1}+q_{L-1}x^{L}+\cdots+q_{1}x^{2}-x+p.$$
Note that $P\in D$, it is easy to see that
$F(x)=0$ has a unique real root in $(0,1)$.
As a consequence $|\lambda I-\hat{M}|=0$ has only one real root larger than 1, so the largest eigenvalue of $\hat{M}$ is larger than one. One can see that 1 is not an eigenvalue of $M$, and the set of eigenvalues of $\hat{M}$ is the union of the set of eigenvalues of ${M}$ and $\{1\}$. So we get $\rho(M)>1$. By Proposition~4.1, $\rho(M)=\lambda_{M}$ is obvious.

Next, for $A\in \mathbb{R}^{L\times L}$, define $\|\cdot\|$ to be the maximum column sum matrix norm, that is: $\|A\|=\max\limits_{1\leq j \leq L}\sum_{i=1}^{L}|a_{ij}|$. Choose an $\varepsilon$ such that both $M-\varepsilon
E$ and $M+\varepsilon E$ are positive matrices. By
Proposition~4.1, both $\lambda_{M}^{\varepsilon-}$ and
$\lambda_{M}^{\varepsilon+}$ are meaningful. By Proposition~4.2, we have
\beqnn
& &d(\sigma(M),\sigma(M-\varepsilon E))\leq4 (\frac{2}{p})^{1-\frac{1}{L}}L^{\frac{1}{L}}\varepsilon^{\frac{1}{L}}\leq C\varepsilon^{\frac{1}{L}},\\
& &d(\sigma(M),\sigma(M+\varepsilon E))\leq4 (\frac{2}{p}+2L)^{1-\frac{1}{L}}L^{\frac{1}{L}}\varepsilon^{\frac{1}{L}}\leq C\varepsilon^{\frac{1}{L}}.
\eeqnn

For such $\varepsilon$, it is obvious that $\lambda_{M}-C\varepsilon^{\frac{1}{L}}>1$. If $\lambda_{M}^{\varepsilon-}<\lambda_{M}-C\varepsilon^{\frac{1}{L}}$, for $\lambda_{M}$, then there exists no permutation such that $d(\sigma(M),\sigma(M-\varepsilon E))\leq C\varepsilon^{\frac{1}{L}}$, because $\lambda_{M}^{\varepsilon-}$ is the largest eigenvalue of $M-\varepsilon E$.

On the other hand, because $M+\varepsilon E>M>0$, we
have $\rho(M+\varepsilon E)\geq\rho(M)>1$ (\cite{hor}, P491, Corollary~ 8.1.19). If
$\lambda_{M}^{\varepsilon+}>\lambda_{M}+C\varepsilon^{\frac{1}{L}}$,
then there exists no permutation such that
$d(\sigma(M),\sigma(M+\varepsilon E))\leq
C\varepsilon^{\frac{1}{L}}$, because $\lambda_{M}$ is the largest
eigenvalue of $M$. \qed

\medskip
\noindent{\it Proof of  Theorem 1.3} ~~ The first
part of the Theorem 1.3 has been proved in Lemma 4.1. Now we
focus on the second part.

{\it Step 1}~ To prove $
\lim_{k\rightarrow \infty}\frac{\log E^{k}T_{k+1}}{k} =\log\lambda_{M}.$  
Let
$$\Phi(i)=\sum_{j=1}^{i}e_{L}M_{i}M_{i-1}\ldots M_{i-j+1}w_{L}',$$
$$\varphi(i)=\sum_{j=1}^{i}M_{i}M_{i-1}\ldots M_{i-j+1}.$$
By formula (\ref{f1.7}),  $E^{i}T_{i+1}=1+\Phi(i).$

Because $P(i)\rightarrow P$, we have $M_{i}\rightarrow M$. By  Lemma~4.1, $\lambda_{M}>1$. For each $\varepsilon$ we defined in Lemma~4.1. $\exists N$, when $n>N$,
$$0<M-\varepsilon E\leq M_{n}\leq M+\varepsilon E.$$
Then we have
\beqnn
&&0<(M-\varepsilon E)^{2}\leq M_{N+2}M_{N+1}\leq (M+\varepsilon E)^{2},\\
&&\cdots\cdots\\
&&0<(M-\varepsilon E)^{k}\leq M_{N+k}M_{N+k-1}\cdots M_{N+1}\leq (M+\varepsilon E)^{k}.
\eeqnn
Summarizing these formulas leads to:
\beqnn
&
&(M-\varepsilon E)^{k}+(M-\varepsilon E)^{k-1}+\cdots +(M-\varepsilon E)\\
&
&\leq M_{N+k}M_{N+k-1}\cdots M_{N+1}+M_{N+k-1}\cdots M_{N+1}+\cdots+M_{N+1}\\
&
&\leq (M+\varepsilon E)^{k}+(M+\varepsilon E)^{k-1}+\cdots +(M+\varepsilon E).
\eeqnn
For the given $N$, Let $A_{N}=M_{N}M_{N-1}\cdots M_{1}+M_{N-1}\cdots M_{1}+\cdots +M_{1}$,
$B_{N}=M_{N}M_{N-1}\cdots M_{1}$.
Then
\beqlb\label{f4.1}
&
&(M-\varepsilon E)^{k}B_{N}+(M-\varepsilon E)^{k-1}B_{N}+\cdots +(M-\varepsilon E)B_{N}+A_{N} \nonumber\\
&
&\leq \varphi(N+k)
\leq (M+\varepsilon E)^{k}B_{N}+(M+\varepsilon E)^{k-1}B_{N}+\cdots +(M+\varepsilon E)B_{N}+A_{N}.
\eeqlb
Now we first consider the ``$\leq$'' part of the
result. From (\ref{f4.1}), we have
\beqnn
\varphi(N+k)\geq (M-\varepsilon E)^{k}B_{N}.
\eeqnn
Then
\beqlb\label{f4.2}
&&\Phi(N+k)\geq (\lambda_{M}^{\varepsilon-})^{k}e_{L}\frac{(M-\varepsilon E)^{k}}{(\lambda_{M}^{\varepsilon-})^{k}}B_{N}w_{L},\nonumber\\
&&\log(\Phi(N+k))\geq k\log\lambda_{M}^{\varepsilon-}+\log e_{L}\frac{(M-\varepsilon E)^{k}}{(\lambda_{M}^{\varepsilon-})^{k}}B_{N}w_{L}.
\eeqlb By {\it Proposition~4.1}, as $k\rightarrow\infty$,
\beqnn
\frac{(M-\varepsilon E)^{k}}{(\lambda_{M}^{\varepsilon-})^{k}}\rightarrow R_{M}^{\varepsilon-}>0,
\eeqnn
as a consequence $e_{L}\frac{(M-\varepsilon E)^{k}}{(\lambda_{M}^{\varepsilon-})^{k}}B_{N}w_{L}$ is bounded in $k$. Thus from (\ref{f4.2}) as  $k\rightarrow\infty$
\beqnn
\varliminf \limits_{k\rightarrow \infty}\frac{\log\Phi(k)}{k}\geq \log\lambda_{M}^{\varepsilon-}.
\eeqnn
Note that from Lemma~4.1, $\lambda_{M}^{\varepsilon-}>\lambda_{M}-C\varepsilon^{\frac{1}{L}}>1$. So
\beqlb\label{f4.8}
\varliminf \limits_{k\rightarrow \infty}\frac{\log\Phi(k)}{k}\geq \log(\lambda_{M}-C\varepsilon^{\frac{1}{L}}).
\eeqlb
For the right ``$\leq$'' part, define $\psi(k)=e_{L}\frac{(M+\varepsilon E)^{k}}{(\lambda_{M}^{\varepsilon+})^{k}}B_{N}w_{L}.$ From (\ref{f4.1}), we have
\beqnn
&&\Phi(N+k)\leq (\lambda_{M}^{\varepsilon+})^{k}(\psi(k)+\psi(k-1)+\cdots+\psi(1)+A_{N}),\\
&&\log \Phi(N+k)\leq k\log(\lambda_{M}^{\varepsilon+})+\log(\psi(k)+\psi(k-1)+\cdots+\psi(1)+A_{N}).
\eeqnn
It's easy to see that $\frac{\log(\psi(k)+\psi(k-1)+\cdots+\psi(1)+A_{N})}{N+k}\rightarrow 0$ as $k\rightarrow \infty$ (because $\lim_{k\to\infty}\psi(k)$ exists). Therefore,
\beqnn
\varlimsup \limits_{k\rightarrow \infty}\frac{\log\Phi(k)}{k}\leq \log\lambda_{M}^{\varepsilon+}.
\eeqnn
By Lemma~4.1, $1<\lambda_{M}^{\epsilon+}<\lambda_{M}+C\epsilon^{\frac{1}{L}}$. So
\beqlb\label{f4.9}
\varlimsup \limits_{k\rightarrow \infty}\frac{\log\Phi(k)}{k}\leq \log(\lambda_{M}+C\varepsilon^{\frac{1}{L}}).
\eeqlb
Combine (\ref{f4.8}) and (\ref{f4.9}) to have
\beqnn
\log(\lambda_{M}-C\varepsilon^{\frac{1}{L}})\leq\varliminf \limits_{k\rightarrow \infty}\frac{\log\Phi(k)}{k}\leq
\varlimsup \limits_{k\rightarrow \infty}\frac{\log\Phi(k)}{k}\leq \log(\lambda_{M}+C\varepsilon^{\frac{1}{L}}).
\eeqnn
Let $\varepsilon\rightarrow 0$, we get ( recall  $E^{i}T_{i+1}=1+\Phi(i)$)
\beqlb\label{t}
\lim_{k\rightarrow \infty}\frac{\log E^{k}T_{k+1}}{k}=\lim_{k\rightarrow \infty}\frac{\log\Phi(k)}{k}=\log\lambda_{M}.
\eeqlb

{\it Step 2}~ It's easy to see that when $P(i)\rightarrow P$ and $P\in D$, the walk is positive recurrence. So the stationary distribution exists and
$\pi (i)=\frac{1}{E^{i}T_{i}}$. By the same method we have used in the proof of Theorem 1.2, we get
\beqlb\label{t2}
E^{i}T_{i}=p(i)E^{i+1}\tau_{i+1}+\sum_{l=1}^{L-1}\vartheta_{l}E^{i-l}T_{i-l+1}+q_{L}(i)E^{i-L}T_{i-L+1}+1
\eeqlb
where $\vartheta_{l}=\sum_{k=l}^{L}q_{k}(i)+p(i)\sum_{j=l}^{L-1}P^{i+1}[(i+1,+\infty),i-j]$. Because the walk is positive recurrence, $E^{i+1}\tau_{i+1}<E^{1}\tau_{1}<\infty$; on the other hand, under our condition, $0<p(i)<1$, $0\leq q_{l}(i)<1$,  $0<\vartheta_{l}<L+1.  $ Then (\ref{t2}) tells us that $E^{i-1}T_{i}$ is the  dominated term for $E^{i}T_{i}$. By (\ref{t}),
\beqnn
\lim_{k\rightarrow \infty}\frac{\log\pi(k)}{k}=-\lim_{k\rightarrow \infty}\frac{\log E^{k}T_{k}}{k}=-\lim_{k\rightarrow \infty}\frac{\log E^{k-1}T_{k}}{k}=-\log\lambda_{M}.
\eeqnn
\qed

\

\noindent{\bf Appendix:} {\it  Analytical proof of Theorem 1.1.}\\

Let $\xi_{0}=0$, $\xi_{1}=1$, and
$\xi_{n}=\sum_{i=1}^{n}\theta_{i}$ for $n\geq 1$, the occupation
time at position 0 before the walk hitting position $n$. By
Theorem~A, we have \beqnn &
&E^{0}\theta_{i}=e_{1}M_{i-1}M_{i-2}\cdots M_{1}u,~~~i\geq 1.
 \eeqnn
\noindent For the system of equations
$\sum_{j=0}^{\infty}P_{ij}y_{j}=y_{i}, i\geq 0$, we will prove
$y_n=E^{0}\xi_{n}$ is the solution, the probabilistic
meaning of which is the expectation of the visiting time at
position 0 before the walk hitting the position $n$. The system of
equations can be rewritten as,
\begin{align*}
y_{1}& =p(1)y_{2}+q_{1}y_{0}, \\
y_{2}& =p(2)y_{3}+q_{1}(2)y_{1}+q_{2}(2)y_{0}, \\
    & \vdots \\
y_{n}& =p(n)y_{n+1}+q_{1}(n)y_{n-1}+q_{2}(n)y_{n-2},\\
& \vdots
\end{align*}
 It is not hard to show that the solution spans a
two dimensional linear space. We can prescribe $y_{0}$ and $y_{1}$
arbitrarily for the initial values, and then all the other $y_{i}$
are determined by these equations. Trivially $y_{i}\equiv 1$ is a
solution. We show that  If $y_n=E^{0}\xi_{n}$,  for $n \geq0$ is
also a solution now.

Setting $x_{0}=0,x_{1}=1$, and for $n\geq2$,
$x_{n}=y_{n}-y_{n-1}$, we obtain:
\begin{align*}
x_{2}&=(\frac{q_{1}(1)}{p_{1}(1)}+\frac{q_{2}(1)}{p_{1}(1)})x_{1}, \nonumber \\
x_{3}&=\frac{q_{2}(2)}{p_{2}(2)}x_{1}+(\frac{q_{1}(2)}{p_{2}(2)}+\frac{q_{2}(2)}{p_{2}(2)})x_{2}, \nonumber\\
& \vdots \nonumber\\
x_{n+1}&=\frac{q_{2}(n)}{p_{2}(n)}x_{n-1}+(\frac{q_{1}(n)}{p_{2}(n)}+\frac{q_{2}(n)}{p_{2}(n)})x_{n}\\
&\vdots \nonumber
\end{align*}
We just need to prove that $x_{n}=E^{0}\theta_{n}$ is
the correspond solution. In probabilistic meaning,  $x_{n}$ is the
expectation of the visiting time at position $0$ in the $n$-th
immigration structure. For the first equation, \beqnn
E^{0}\theta_{2}=e_{1}
M_{1}u=\frac{q_{1}(1)}{p_{1}(1)}+\frac{q_{2}(1)}{p_{1}(1)}=(\frac{q_{1}(1)}{p_{1}(1)}+\frac{q_{2}(1)}{p_{1}(1)})E^{0}\theta_{1}.
\eeqnn \noindent For the $n$-th equation, let $K_{n}=M_{n-2}
M_{n-3}\cdots
M_{1}u,b_{n}(1)=\frac{q_{1}(n)}{p(n)},b_{n}(2)=\frac{q_{2}(n)}{p(n)}.$
\beqnn E^{0}\theta_{n+1}
&=&e_{1}M_{n}M_{n-1}K_{n}\\
&=&\Big(b_{n}(1)b_{n-1}(1)+b_{n}(2)b_{n-1}(2)+b_{n}(2),~b_{n}(1)b_{n-1}(2)+b_{n}(2)b_{n-1}(2)\Big)K_{n}.
\eeqnn
 On the other hand, \beqnn
& &\frac{q_{2}(n)}{p_{2}(n)}E^{0}\theta_{n-1}+(\frac{q_{1}(n)}{p_{2}(n)}+\frac{q_{2}(n)}{p_{2}(n)})E^{0}\theta_{n}\\
&=&(b_{n}(2)e_{1}+(b_{n}(1)+b_{n}(1))e_{1}M_{n-1})K_{n}\\
&=&((b_{n}(2),0)+(b_{n}(1)+b_{n}(2),0)M_{n-1})K_{n}\\
&=&\Big((b_{n}(2),0)+(b_{n}(1)+b_{n}(2),0)\left(
                                       \begin{array}{cc}
                                         b_{n}(1) & b_{n}(2) \\
                                         1+b_{n}(1) & b_{n}(2) \\
                                       \end{array}
                                    \right)\Big)K_{n}\\
&=&\Big(b_{n}(1)b_{n-1}(1)+b_{n}(2)b_{n-1}(2)+b_{n}(2),b_{n}(1)b_{n-1}(2)+b_{n}(2)b_{n-1}(2)\Big)K_{n},
\eeqnn
 which means that $x_{n}=E^{0}\theta_{n}$ is the correspond solution of the equations above.
So $y_n=E^{0}\xi_{n}$ is the solution of
$\sum_{j=0}^{\infty}P_{ij}y_{j}=y_{i}, i\geq 0$. It is evidence
that
  $y_{n}\equiv 1$ are another solutions of the system. Therefore, the general solution is given
  by
\beqnn
(y_{0},y_{1},\ldots,y_{n}\ldots)=\alpha\cdot(1,1,\ldots,1,\ldots)'+\beta\cdot(E^{0}\xi_{0},E^{0}\xi_{1},\ldots,E^{0}\xi_{n}\ldots)',
\eeqnn
where $\alpha,\beta,\in \mathbb{R}$, and a nonconstant bounded solution exists if and only if $E^{0}\xi_{n}$
is bounded, i.e., $\sum_{k=1}^{\infty}e_{1}M_{k}M_{k-1}\cdots M_{1}u< \infty$ \qed\\

\noindent{\bf Remark}  The probability meaning of $\theta_{i}$,
$i>0$, is the expectation of the local time at zero of
the walk starting at $i-1$ before $T_{i}$, which is the
motivation of the proof.

\end{document}